\documentclass[12pt]{amsart}

\newtheorem{lem}[equation]{Lemma}
\newtheorem{thm}[equation]{Theorem}
\newtheorem{prop}[equation]{Proposition}
\newtheorem{defn}[equation]{Definition}

\newtheorem{cor}[equation]{Corollary}
\newtheorem{rem}[equation]{Remark}
\usepackage[draft]{graphicx}
\usepackage{amssymb, amsmath}
\usepackage[english]{babel}
\usepackage{amsfonts}
\usepackage{color}

\newcommand{\dw}{\downarrow}
\newcommand{\uw}{\uparrow}

\newcommand{\ch}{\chi}
\newcommand{\mtr}{\mathrm}

\newcommand{\ncm}{\newcommand}
\ncm{\np}{\newpage}
\ncm{\ebl}{\end{thebibliography}}
\ncm{\bbl}{\begin{thebibliography}}
\ncm{\chd}{_{ _{\ch}}}
\ncm{\ald}{_{ _{\al}}}
\ncm{\cP}{\mathcal{P}}
\ncm{\ei}{e_i}
\ncm{\eij}{e_{i,\;j}}
\ncm{\bt}{\begin{thm}}
\ncm{\bdef}{\begin{defn}}
\ncm{\edf}{\end{defn}}
\ncm{\et}{\end{thm}}
\ncm{\bc}{\begin{cor}}
\ncm{\bl}{\begin{lem}}
\ncm{\el}{\end{lem}}
\ncm{\bpf}{\begin{proof}}
\ncm{\epf}{\end{proof}}
\ncm{\ec}{\end{cor}}
\ncm{\er}{\end{rem}}
\ncm{\br}{\begin{rem}}
\ncm{\bn}{\begin}
\ncm{\bp}{\begin{prop}}
\ncm{\ep}{\end{prop}}
\ncm{\bd}{\begin{document}}
\ncm{\ed}{\end{document}}
\ncm{\beq}{\begin{equation}}
\ncm{\beqn}{\begin{equation*}}
\ncm{\eeq}{\end{equation}}
\ncm{\eeqn}{\end{equation*}}
\ncm{\bea}{\begin{eqnarray}}
\ncm{\eea}{\end{eqnarray}}
\ncm{\beanon}{\begin{eqnarray*}}
\ncm{\eeanon}{\end{eqnarray*}}\ncm{\ek}{\eps|_K}\ncm{\diez}{\#}
\ncm{\bwt}{\bowtie}
\ncm{\cC}{\mtc{C}}
\ncm{\cX}{\mtc{X}}
\ncm{\wt}{\widetilde}
\ncm{\sg}{\sigma}\ncm{\Rep}{\mathrm{Rep}}
\ncm{\Irr}{\mathrm{Irr}}\ncm{\X}{\mathcal{X}}
\ncm{\cA}{\mathcal{A}}
\ncm{\HKer}{\mtr{HKer}}
\ncm{\LKER}{\mtr{LKER}}
\ncm{\aad}{\mtr{ad}}
\ncm{\Dr}{\mtr{D}}
\ncm{\cD}{\mathcal{D}}
\ncm{\G}{\mathcal{G}}
\ncm{\Dc}{\mtc{D}}
\ncm{\E}{\mtc{E}}
\ncm{\fp}{\mtr{FPdim}}
\ncm{\Vc}{\mtr{Vec}}
\ncm{\cK}{\mtc{K}}
\ncm{\cM}{\mtc{M}}
\ncm{\cE}{\mtc{E}}
\ncm{\cS}{\mtc{S}}
\ncm{\End}{\mtr{End}}

\newcommand{\DOT}{\setlength{\unitlength}{1pt}\begin{picture}(2.5,2)
                  (1,1)\put(2,3.5){\circle*{2}}\end{picture}}

\newcommand{\coh}{{\rm H}}
\newcommand{\HH}{{\rm HH}}
\newcommand{\ra}{\rightarrow}
\newcommand{\codim}{{\rm codim}}
\newcommand{\Wedge}{\textstyle\bigwedge}
\newcommand{\C}{\mathbb C}
\newcommand{\R}{\mathbb R}
\newcommand{\Z}{\mathbb Z}
\newcommand{\M}{\mathcal M}
\newcommand{\ot}{\otimes}
\newcommand{\mtc}{\mathcal}
\newcommand{\lam}{\lambda}
\newcommand{\lb}{\label}
\newcommand{\Lam}{\Lambda}
\newcommand{\lbd}{\Lambda}
\newcommand{\sig}{\sigma}
\newcommand{\al}{\alpha}
\newcommand{\eps}{\epsilon}
\newcommand{\en}{\end}
\newcommand{\teta}{\theta}
\newcommand{\nhs}{normal Hopf subalgebra}
\newcommand{\ul}{\underline}
\newcommand{\mr}{\mathrm}
\newcommand{\rh}{\rightharpoonup}
\newcommand{\lh}{\leftharpoonup}
\newcommand{\sub}{\subsection}
\newcommand{\mc}{\mathcal}
\newcommand{\D}{\Delta}
\ncm{\rep}{\mtr{Rep}}
\ncm{\btw}{\bowtie}
\ncm{\cd}{\mtc{D}}
\ncm{\cs}{\mtc{S}}
\ncm{\bq}{\beq}
\ncm{\eq}{\eeq}
\ncm{\cop}{\mtr{cop}}
\newcommand{\nc}{\newcommand}
\newcommand{\gm}{\gamma}
\newcommand{\beqarn}{\begin{eqnarray*}}
\newcommand{\eeqarn}{\end{eqnarray*}}
\ncm{\cc}{\mtc{C}}
\ncm{\blue}{\textcolor[rgb]{.00, .00, 1.00}}
\ncm{\red}{\textcolor[rgb]{1.00, .00, .00}}
\ncm{\md}{\medbreak}
\ncm{\mi}{\mtr{I}}
\ncm{\cZ}{\mtc{Z}}\ncm{\xra}{\xrightarrow}
\ncm{\cb}{\mtc{B}}\ncm{\ca}{\mtc{A}}
\ncm{\cm}{\mathcal{M}}
\ncm{\co}{\mtc{O}}
\ncm{\bne}{\begin{enumerate}}
\ncm{\ene}{\end{enumerate}}
\ncm{\os}{\oplus}
\ncm{\stab}{\mtr{Stab}}
\ncm{\onh}{On the other hand}
\title[Semisimple Hopf algebras]
{ New examples of Green functors arising from representation theory of semisimple Hopf algebras}
\thanks{This work was supported by a grant of the Romanian National Authority for Scientific Research, CNCS Ð UEFISCDI, project number PN-II-RU-TE-2012-3-0168.} 
\ncm{\inv}{^{-1}}
\author{Sebastian  Burciu}

\begin{document}
\maketitle
\begin{abstract}
A general Mackey type decomposition for representations of semisimple Hopf algebras  is investigated. We show that such a decomposition occurs in the case that the module is induced from an arbitrary Hopf subalgebra and it is restricted back to a group subalgebra. Some other examples when such a decomposition occurs are also constructed. They arise from gradings on the category of corepresentations of a  semisimple Hopf algebra and provide new examples of Green functors  in the literature.
\end{abstract}

\section{Introduction and Main Results}\lb{imr}
Mackey's decomposition theorem of induced modules from subgroups is a very important tool in the representations theory of finite groups. This decomposition describes the process of an induction composed with a restriction in terms of the reverse processes consisting of restrictions followed by inductions. More precisely, if $G$ is a finite group, $M$ and $N$ two subgroups of $G$  and $V$ a finite dimensional $k$-linear representation of $M$ then the well known Mackey's decomposition states that there is an isomorphism of $kN$-modules:
\beq
V\uw^{kG}_{kM}\dw^{kG}_{kN}\;
\xra{\delta_V}
\bigoplus_{x \in M\backslash G\slash N}
k[ N]\ot _{k[\;^xM\cap N]}\;^xV.
\eeq
Here $\;^xM:=xMx^{-1}$ is the conjugate subgroup and $\;^xV:=V$ is the conjugate $\;^xM$-representation defined by $(xmx^{-1}).v:=m.v$ for all $m \in M$ and $v \in V$. The direct sum is indexed by a set of representative group elements of $G$ for all double cosets $M\backslash G\slash N$ of $G$ relative to the two  subgroups $M$ and $N$. Note that the inverse isomorphism of $\delta_V$ is given on each direct summand by the left multiplication operator $n\ot_{kN\cap k\;^xM} v\mapsto nx \ot_{kM}v$, see \cite[Proposition 22]{Se76}. 
\md
The goal of this paper is to investigate a similar Mackey type decomposition for the induced modules from Hopf subalgebras of semisimple Hopf algebras and restricted back to other Hopf subalgebras. In order to do this, we use the corresponding notion of a double coset relative to a pair of Hopf subalgebras of a semisimple Hopf algebra that was introduced by the author in \cite{coset}. We also have to define a conjugate Hopf subalgebra corresponding to the notion of a conjugate subgroup. For any Hopf subalgebra $K\subseteq H$ of a semisimple Hopf algebra $H$ and any simple subcoalgebra $C$ of $H$ we define the conjugate Hopf subalgebra $^CK$ of $K$ in Proposition \ref{conjdef}. This notion corresponds to the notion of conjugate subgroup from the above decomposition.  In order to deduce that $^CK$ is a Hopf subalgebra of $H$ we use several crucial results from \cite{NZ} concerning the product of two subcoalgebras of a semisimple Hopf algebra as well as Frobenius-Perron theory for nonnegative matrices. \md
Using these tools we can prove one of the following main results of this paper: 
\bt\lb{grpcase}
Let $K\subseteq H$ be a Hopf subalgebra of a semisimple Hopf algebra and $M$ a finite dimensional $K$-module. Then for any subgroup $G\subseteq G(H)$ one has a canonical isomorphism of $kG$-modules
\beq\lb{stg}
M\uw^H_K\dw^H_{kG}\; \;\xra{\delta_M} \bigoplus_{C \in kG\backslash H/K}(kG\ot_{
\;kG_C}\;^CM).
\eeq
Here $G(H)$ is the group of grouplike elements of $H$ and the subgroup $G_C\subseteq G$ is determined by $kG\;\cap\; ^CK=kG_C$. The conjugate module $\:^CM$ is defined by $\:^CM:=CK\ot_KM$.
\et
\ncm{\sem}{semisimple\;}
As in the classical group case the homomorphism $\delta_M$ is the inverse of a natural homomorphism $\pi_M$ which is constructed by the left multiplication on each direct summand. It is not difficult to check (see Theorem \ref{mack} below) that in general, for any two Hopf subalgebras $K, L\subseteq H$ the left multiplication homomorphism $\pi_M$ is always an epimorphism:
\bt\label{mack}
Let $K$ and $L$ be two Hopf subalgebras of a semisimple Hopf algebra $H$. For any finite dimensional left $K$-module $M$ there is a canonical epimorphism of $L$-modules
\beq\lb{st}
\oplus_{C \in L\backslash H/K}(L\ot_{L\cap 
\;^CK}\;^CM) \xra{\pi_M} M\uw^H_K\dw^H_L
\eeq
given on components  by $l \ot_{L\cap 
\;^CK} v\mapsto lv$ for any $l \in L$ and any $v \in \;^CM$.
Here  the conjugate module $^CM$ is defined as above by $^CK:=CK\ot_KM$.
\et
\md
We remark that there is  a similar direction in the literature in the paper \cite{fes}. In this paper the author considers a similar decomposition but for pointed Hopf algebras instead of semisimple Hopf algebras. Also, in \cite{lin} the author proves a similar result for some special  Hopf subalgebras of quantum groups at roots of $1$.
\md Another particular situation of Mackey's decomposition can be found in \cite{coset}. In this paper it is proven that for pairs of Hopf subalgebras that generate just one double coset subcoalgebra, the above epimorphism $\pi_M$ from Theorem \ref{mack} is in fact an isomorphism. In both papers, the above homomorphism $\pi_M$ is given by left multiplication.
\md
\bn{defn}
We say that $(L,K)$ is a Mackey pair of Hopf subalgebras  of $H$ if the above left mutliplication homomorphism $\pi_M$ from Theorem \ref{mack} is an isomorphism for any finite dimensional left $K$-module $M$.
\end{defn}
Then Theorem \ref{grpcase} above states that $(kG, K)$ is a Mackey pair for any Hopf subalgebra $K\subset H$ and any subgroup $G\subset G(H)$.
 Moreover  in Theorem \ref{norm} it is shown that for any normal Hopf subalgebra $K$ of $H$ the pair $(K,K)$ is a Mackey pair.  This allows us to prove a new formula (see Proposition \ref{irredind}) for the restriction of an induced module from a normal Hopf subalgebra which substantially improves Proposition 5.12 of \cite{coset}. It also gives a criterion for an induced module from a normal Hopf subalgebra to be irreducible generalizing a well known criterion for group representations, see for example \cite[Corollary 7.1]{Se76}.
\md
For any semisimple Hopf algebra $H$, using the universal grading of the fusion category $\rep(H^*)$ we construct in Section \ref{newgf} new Mackey pairs of Hopf subalgebras of $H$. In turn, this allows us to define a Green functor on the universal group $G$ of the category of representations of $H^*$. For $H=kG$ one obtains in this way the usual Green functor \cite{green}. As in group theory, this new Green functor can be used to determine new properties of the Grothendieck ring of a semisimple Hopf algebra.
\md In the last section we prove the following tensor product formula for two induced modules from a Mackey pair of Hopf subalgebras.
\bt \lb{tp}Suppose that $(L,K)$ is a Mackey pair of Hopf subalgebras of a semisimple Hopf algebra $H$. Then for any $K$-module $M$ and any $L$-module $N$ one has a canonical isomorphism:
\beq
M\uw_K^H\ot N\uw_L^H\;\; \xrightarrow{\cong} \bigoplus_{C \in L\backslash H/K}((CK\ot_KM)\dw^{\;^CK}_{\;L\cap \;^CK}\ot N\dw^{\;L}_{\;L\cap \;^CK})\uw^{\;H}_{\;L\cap \;^CK}
\eeq
\et
This generalizes a well known formula for the tensor product of two induced  group representations given for example in \cite{benson}.
\md
This paper is structured as follows. In the first section we recall the basic results on coset decomposition for Hopf algebras. The second section contains the construction for the conjugate Hopf subalgebra generalizing the conjugate subgroup of a finite group. These results are inspired from the treatment given in \cite{d3k}. A general characterization for the conjugate Hopf subalgebra is given in Theorem \ref{larg}. This theorem is automatically satisfied in the group case.
In the third section we prove Theorem \ref{mack}.  In the next section we prove Theorem \ref{grpcase}.
We also show that for any semisimple Hopf algebra there are some canonical associated Mackey pairs arising from the universal grading of the category of finite dimensional corepresentations (see Theorem \ref{grd}). Necessary and sufficient conditions for a given pair to be a Mackey pair are given in terms of the dimensions of the two Hopf subalgebras of the pair and their conjugate Hopf subalgebras. In Section \ref{nrm} we prove that for a normal Hopf subalgebra $K$ the pair $(K,K)$ is always a Mackey pair.
In the last subsection \ref{tpf} we prove the tensor product formula from Theorem \ref{tp}.
\md
We work over an algebraically closed field $k$ of characteristic zero. We use Sweedler's notation $\Delta$  for comultiplication but with the sigma symbol dropped. All the other Hopf algebra notations of this paper are the standard ones, used for example in \cite{montg}.
\section{ Double coset decomposition for Hopf subalgebras of semisimple Hopf algebras}\lb{dbc}
\sub{Conventions:} 
Throughout this paper $H$ will be a semisimple Hopf algebra over $k$ and $\Lam_H \in H$ denotes its idempotent integral ($\eps(\Lam_H)=1$). It follows that $H$ is also cosemisimple \cite{Lard}. If $K$ is a Hopf
subalgebra of $H$ then $K$ is also a semisimple and cosemisimple Hopf algebra \cite{montg}. For any two subcoalgebras $C$ and $D$ of $H$ we denote by $CD$ the subcoalgebra of $H$ generated as a $k$-vector space by all elements of the type $cd$ with $c \in C$ and $d \in D$.
\md
Let $G_0(H)$ be the Grothendieck group of the category of left $H$-modules. Then since $H$ is a Hopf algebra the group $G_0(H)$ has a ring structure under the tensor product of modules. Then the character ring $C(H):=G_0(H)\ot_{\mathbb{Z}}k$ is a semisimple subalgebra of  $H^*$ \cite{Zhc}. Denote by $\Irr(H)$ the set of all irreducible characters of $H$. 
Then $C(H)$ has a basis consisting of the irreducible characters $\ch \in \Irr(H)$. Also $C(H)$ coincides to the space $\mtr{Cocom}(H^*)$ of cocommutative elements of $H^*$. By duality, the character ring $C(H^*)$ of the dual Hopf algebra $H^*$ is a semisimple Hopf subalgebra of $H$ and $C(H^*)=\mtr{Cocom}(H)$. If $M$ is a finite dimensional $H$-module with character $\ch$ then the linear dual $M^*$ becomes a left $H$-module with character $\ch^*:=\ch\circ S$.
\sub{The subcoalgebra associated to a comodule}\lb{subc}
Let $W$ be any right $H$-comodule. Since $H$ is finite dimensional it follows that $W$ is a left $H^*$-module via the module structure $f.w=f(w_1)w_0$,  where $\rho (w)=w_0\ot w_1$ is the given right $H$-comodule structure of $W$. Then one can associate to $W$ the coefficient subcoalgebra denoted by $C_W$ \cite{Lar}. Recall that $C_W$ is the minimal subcoalgebra $C$ of $H$ with the property that $\rho(W)\subset W\ot C$. Moreover, it can be shown that $C_W=(\mtr{Ann}_{H^*}(W))^{\perp}$ and $C_W$ is called the subcoalgebra of $H$ associated to the right $H$-comodule $W$.  
If $W$ is a simple right $H$-comodule (or equivalently $W$ is an irreducible $H^*$-module) then the associated subcoalgebra $C_W$ is a co-matrix coalgebra. More precisely, if $\dim W=q$ then $\dim C_W=q^2$ and it has a a $k$-linear basis given by $x_{ij}$ with $1\leq i,j\leq q$ . The coalgebra structure of $C_W$ is then given  by $\D(x_{ij})=\sum_{l}x_{il}\ot x_{lj}$ for all $1\leq i, j \leq q$. Moreover the irreducible character $d \in C(H^*)$ of $W$ is given by formula $d=\sum_{i=0}^qx_{ii}$. It is easy to check that $W$ is an irreducible $H^*$-module if and only if $C_W$ is a simple subcoalgebra of $H$. This establishes a canonical bijection between the set $\Irr(H^*)$ of simple right $H^*$-comodules and the set of simple subcoalgebras of $H$. For any irreducible character $d \in \Irr(H^*)$ we also use the notation $C_d$ for the simple subcoalgebra of $H$ associated to the character $d$ (see \cite{Lar}).

Recall also that if $M$ and $N$ are two right $H$-comodules then $M\ot N$ is also a comodule with $\rho(m\ot n)=m_0\ot n_0 \ot m_1n_1$.

\br\lb{princ}
For a simple subcoalgebra $C\subset H$ we denote by $M_C$ the simple $H$-comodule associated to $C$. Following \cite{NZ}, if $C$ and $D$ are simple subcoalgebras of a semisimple Hopf algebra $H$ then the simple comodules entering in the decomposition of $M_C\ot M_D$ are in bijection with the set of all simple subcoalgebras of the product  subcoalgebra $CD$ of $H$. Moreover, this bijection is given by $W\mapsto C_W$ for any simple subcomodule $W$ of $M_C\ot M_D$.
\er
\subsection{Double coset decomposition for Hopf subalgebras}
In this subsection we recall the basic facts on double cosets of semisimple Hopf algebras developed in \cite{coset}. Let $L$ and $K$ be two Hopf subalgebras of $H$.  As in \cite{coset} one can define an equivalence
relation $r_{ _{L,\;K}}^H$ on the set of simple subcoalgebras of $H$ as following: $C
\sim D$ if  $C \subset \mathrm{LDK}$. The fact that $r_{ _{L,\;K}}^H$ is an equivalence relation is proven in \cite{coset}. In this paper it is shown that $C\sim D$ if and only if $LDK=LCK$ as subcoalgebras of $H$. We also have the following:
\bn{prop}\lb{dc}
If $C$ and $D$ are two simple subcoalgebras of $H$ then the following are equivalent:

1) $C\sim\; D$

2) $LCK=LDK$

3) $\Lam_LC\Lam_K=\Lam_LD\Lam_K$
\end{prop}
\bn{proof} First assertion is equivalent to the second by Corollary 2.5 from \cite{coset}. Clearly $(2)\Rightarrow (3)$ by left multiplication with $\Lam_K$ and right multiplication with $\Lam_L$. It will be shown that  $(3)\Rightarrow (1)$. One has the following decomposition:
\begin{equation*}
H=\bigoplus_{i=1}^lLC_iK
\end{equation*}
where $C_1,\;\cdots,\;C_l$ are representative subcoalgebras for each equivalence class of $r_{ _{K,\;L}}^H$.

It follows that $\Lambda_LH\Lambda_K=\bigoplus_{i=1}^l\Lambda_LC_i\Lambda_K$. Thus if $C \nsim D$ then $\Lam_LC\Lam_K \cap \Lam_LD\Lam_K=0$ which proves $(1)$.
\end{proof}
\bn{rem}\label{nf} 
The above Proposition shows that for any two simple subcoalgebras $C$ and $D$ of $H$ then either $LCK=LDK$ or $LCK \cap LDK=0$. Therefore for any subcoalgebra $D \subset LCK$ one has that $LCK=LDK$. In particular, for $L=k$, the trivial Hopf subalgebra, one has that  $D \subset CK$ if and only if $DK=CK$.
\end{rem}
\subsubsection{Notations}\lb{not} For the rest of the paper we denote by $L\backslash H\slash K$ the set of double cosets $LCK$ of $H$ with respect to $L$ and $K$. Thus the elements $LCK$ of $L\backslash H\slash K$ are given by a choice of representative of simple subcoalgebras in each equivalence class of $r_{ _{L,\;K}}^H$. Similarly, we denote by $H/K$ be the set of right cosets $CK$ of $H$ with respect to $K$. This corresponds to a choice of a representative simple subcoalgebra in each equivalence class of $r_{ _{k,\;K}}^H$.
\br
As noticed in \cite{coset} one has that $LCK \in \cm^H_K$ and therefore $LCK$ is a free right $K$-module. Similarly $LCK \in \;^H_L\cm$ and therefore $LCK$ is also a free left $L$-module.
\er
By Corollary 2.6 of \cite{coset} it follows that two simple subcoalgebras $C$ and $D$ are in the same double coset of $H$ with respect to $L$ and $K$ if and only if 
\beq\lb{twisted}
\Lam_L\frac{c}{\eps(c)}\Lam_K=\Lam_L\frac{d}{\eps(d)}\Lam_K.
\eeq
where $c$ and $d$ are the irreducible characters of $H^*$ associated to the simple subcoalgebras $C$ and $D$. In particular for $L=k$, the trivial Hopf subalgebra, it follows that $CK=DK$ if and only if
\beq\lb{onesided}
c\Lam_K=\frac{\eps(c)}{\eps(d)}d\Lam_K.
\eeq
\subsection{Principal eigenspace for $<C>$}
For a simple subcoalgebra $C$ we denote by $<C>$ the Hopf subalgebra of $H$ generated by $C$. If $d$ is the character associated to $C$ we also denote this Hopf subalgebra by $<d>$.
\md
\subsubsection{Frobenius-Perron theory for nonnegative matrices}
Next we will use the Frobenius-Perron theorem for matrices with
nonnegative entries (see \cite{F}). If $A\geq 0$ is such a matrix then $A$
has a positive eigenvalue $\lambda$ which has the biggest absolute
value among all the other eigenvalues of $A$. The eigenspace
corresponding to $\lambda$ has a unique vector with all entries
positive. $\lambda$ is called the principal eigenvalue of $A$ and the
corresponding positive vector is called the principal vector of $A$. Also the eigenspace of $A$ corresponding to $\lam$ is called the principal eigenspace of  the matrix $A$. \md
For an irreducible character $d \in \Irr(H^*)$ let $L_d$ be the linear operator on $C(H^*)$ given by left multiplication by $d$.
Recall \cite{coset} that $\eps(d)$ is the Frobenius-Perron eigenvalue of the nonnegative matrix associated to the operator $L_d$ with respect to the basis given by the irreducible characters of $H^*$. In analogy with Frobenius-Perron theory, for a subcoalgebra $C$ with associated character $d$ we call the space of eigenvectors of $L_{d}$ corresponding to the eigenvalue $\eps(d)$ as the principal eigenspace for $L_{d}$.

Next Corollary is a particular case of Theorem 2.4 of \cite{coset}.

\bc\label{print}The principal eigenspace of $L_{\Lam_K}$ is $\Lam_KC(H^*)$ and it has a $k$-linear basis given by $\Lam_Kd$ where $d$ are the characters of a set of representative simple coalgebras for the right cosets of $K$ inside $H$.
\ec
Using this we can prove the following:
\bt\label{eigenv}Let $C$ be a subcoalgebra of a semisimple Hopf algebra $H$ with associated character $d \in C(H^*)$. Then the principal eigenspaces of  $L_d$ and $L_{ _{\Lam_{<d>}}}$ coincide.
\et
\begin{proof}
Let $V$ be the principal eigenspace of $L_{ _{\Lam_{<d>}}}$ and $W$ be the principal eigenspace of $L_d$. Then by Corollary \ref{print} one has that $V=\Lam_{<d>}C(H^*)$. Since $d\Lam_{<d>}=\eps(d)\Lam_{<d>}$ then clearly $V\subseteq W$. On the other hand since $\Lam_{<d>}$ is a polynomial with rational coefficients in $d$ (see Corollary 19 of \cite{nri}) it also follows that $W\subseteq V$.
\end{proof}

\subsection{Rank of cosets} Let $K$ be a Hopf subalgebra of a semisimple Hopf algebra $H$. Consider the equivalence relation $r_{k,\;K}^H$ on the set  $\Irr(H^*)$ of  simple subcoalgebras of $H$. As above one has  $C\sim D$ if and only if $CK=DK$.
Therefore 
\beq\lb{decomp}
H=\oplus_{C \in H/K }CK.
\eeq

\bl\lb{onek}
The equivalence class under $r^H_{k,K}$ of the trivial subcoalgebra $k$ is the set of all simple subcoalgebras of $K$.
\el
\begin{proof}
Indeed suppose that $C$ is a simple subcoalgebra of $H$ equivalent to the trivial subcoalgebra $k$. Then $CK=kK=K$ by Proposition \ref{dc}. Therefore $C \subset CK=K$. Conversely, if $C\subset K$ then $CK\subset K$ and since $CK\in \mtc{M}^H_K$ it follows that $CK=K$. Thus $C\sim k$.
\end{proof}

\bn{prop}\lb{frns}
If $D$ is a simple subcoalgebra of a semisimple Hopf algebra $H$ and $e \in K$ is an idempotent then $$DK\ot_KKe \cong DKe$$ as vector spaces.
\end{prop}

\bn{proof}
Since $H$ is free right $K$-module one has that the map $$\phi : H\ot_K Ke \ra He, \;\;\; h\ot_Kre \mapsto hre$$ is an isomorphism of $H$-modules. Using the above decomposition \eqref{decomp} of $H$ and the fact that $DK$ is a free right $K$-module note that $\phi$ sends $DK\ot_KKe$ to $DKe$.
\end{proof}

\bn{cor} Let $K$ be a Hopf subalgebra of a semisimple Hopf algebra $H$. For any simple subcoalgebra $C$ of $H$ one has that the rank of $CK$ as right $K$-module is $\dim_k C\Lam_K$.  
\end{cor}

\begin{proof}
Put $e=\Lam_K$ the idempotent integral of $K$ in the above Proposition.
\end{proof}

\subsection{Frobenius-Perron eigenvectors for cosets}\lb{eignv} 
 Let $T$ be the linear operator given by right multiplication with $\Lam_K$ on the character ring $C(H^*)$. 
 
 \br\lb{coseignv}
 Using Proposition 2.5 from \cite{coset} it follows that the largest (in absolute value) eigenvalue of $T$ equals $\dim K$. Moreover a basis of eigenvectors corresponding to this eigenvalue is given by $c\Lam_K$ where the character $c \in \Irr(H^*)$ runs through  a set of irreducible characters  representative for all the right cosets $CK\in H/K$.
 \er

\section{The conjugate Hopf subalgebra $\;^CK$}\lb{conj}
Let as above $K$ be a Hopf subalgebra of a semisimple Hopf algebra $H$.
 For any simple subcoalgebra $C$ of $H$ in this section we construct the conjugate Hopf subalgebra $\;^CK$ appearing in Theorem \ref{st}.   If $c \in \Irr(H^*)$ is the associated irreducible character  of $C$ then consider the following subset of $\Irr(H^*)$:
\beq\lb{chcj}
^cK=\{d\in \mtr{Irr}(H^*)\;|\;dc\Lam_K=\eps(d)c\Lam_K\;\}
\eeq
where as above $\Lam_K\in K$ is the idempotent integral of $K$.
\md
Recall from \cite{NZ} that a subset $X \subset \mtr{Irr}(H^*)$ is closed under multiplication if for every two elements $c, d \in X$ in the decomposition of the product $cd=\sum_{e \in \mtr{Irr}(H^*)}m_{c,d}^{e}e$ one has $e\in X$ whenever $m_e\neq 0$. Also a subset $X \subset \mtr{Irr}(H^*)$ is closed under $``\;^*\;"$ if $x^* \in X$ for all $x \in X$. 
\md
Following \cite{NZ} any subset $X \subset \Irr(H^*)$ closed under multiplication generates a subbialgebra $H(X)$ of $H$ defined by 
\beq
H(X):=\oplus_{x \in X}C_x.
\eeq

Moreover if the set $X$ is also closed under  $"\;^*\;"$ then $H(X)$ is a Hopf subalgebra of $H$.
\br\lb{subbualg}Since in our case $H$ is finite dimensional, it is well known that any subbialgebra of $H$ is also a Hopf subalgebra. Therefore in this case any set $X$ of irreducible characters closed under product is also closed under $"^*"$.
\er
\bn{prop}\lb{conjdef}
The set $\;^cK \subset \Irr(H^*)$ is closed under multiplication and $"^*"$ and it generates a Hopf subalgebra $^CK$ of $H$. Thus
\beq
^CK=\oplus_{d \in ^cK}C_d
\eeq
\end{prop}
\bn{proof}
Suppose that $D$ and $D'$ are two simple subcoalgebras of $H$ whose irreducible characters satisfy $d, d' \in \;^cK$. Then one has $dd'c\Lam_K=\eps(dd')c\Lam_K$. On the other hand suppose that 
\beq
dd'=\sum_{e \in \Irr(H^*)}m_{d,d'}^ee.
\eeq
 Then $\eps(dd')c\Lam_K=dd'c\Lam_K=\sum_{e \in \Irr(H^*)}m_{d,d'}^eec\Lam_K$
and Remark \ref{coseignv} implies that $ec\Lam_K$ is a scalar multiple of $c\Lam_K$ for any $e $ with $m_{d,d'}^e\neq 0$. Therefore $ec\Lam_K=\eps(e)c\Lam_K$ and $e \in \;^cK$. This shows that $\;^CK$ is a subbialgebra of $H$ and by Remark \ref{subbualg} a Hopf subalgebra of $H$. 
\end{proof}
Sometimes the notation $^CK$ will also be used for $^cK$ where $c \in \Irr(H^*)$ is the irreducible character associated to the simple subcoalgebra $C$.
\md\noindent
The notion of conjugate Hopf subalgebra $^CK$ is motivated by the following Proposition:
\bp Let $H$ be a semisimple Hopf algebra over $k$. If the simple subcoalgebra $C$ is of the form $C=kg$ with $g \in G(H)$ a group-like element of $H$ then $\;{^{C}}K=gKg^{-1}$.
\ep

\begin{proof}
Indeed, suppose that $D \in \;^CK$.  If $d$ is the associated irreducible character of $D$ then by definition it follows that $dg\Lam_K=g\Lam_K$. Thus $g^{-1}dg\Lam_K=\Lam_K$. Therefore the simple subcoalgebra $g^{-1}Dg$ of $H$ is equivalent to the trivial subcoalgebra $k$. Then using Lemma \ref{onek} one has that $g^{-1}Dg \subset K$  and therefore $\;^CK \subset gKg^{-1}$. The other inclusion $gKg^{-1} \subset \;^CK$ is obvious.\end{proof}
\br
In particular for $H=kG$ one has that $\;^Ck[M]=k[\;^xM]$ where $x \in G$ is given by $C=kx$.
\er
\br\lb{comm}
\bne
\item Using Remark \ref{princ} it follows from the definition of conjugate Hopf subalgebra that $CK$ is always a left $\;^CK$-module.
\item Note that if $C(H^*)$ is commutative then $\;^CK\supseteq K$. Indeed for any $d \in \Irr(K^*)$ one has $d\Lam_K=\eps(d)\Lam_K$ and therefore $dc\Lam_K=cd\Lam_K=\eps(d)c\Lam_K$.
\item If $K$ is a normal Hopf subalgebra of $H$ then since $\Lam_K$ is a central element in $H$ by same argument it also follows that $\;^CK\supseteq K$.
\ene
\er
\subsection{Some properties of the conjugate Hopf subalgebra}
%
\bp\lb{modstr} Let $H$ be a semisimple Hopf algebra and $K$ be a Hopf subalgebra of $H$. Then for any simple subcoalgebra $C$ of $H$ one has that $\;^CK$ coincides to the maximal Hopf subalgebra $L$ of $H$ with the property $LCK=CK$.
\ep
\begin{proof}
The equality $\;^CKCK=CK$ follows from the character equality $\Lam_{\;^CK}c\Lam_K=\eps(\Lam_{\;^CK})c\Lam_K$ and Reamark \ref{princ}. Conversely, if $LCK=CK$ by passing to the regular $H^*$-characters and using Equation \eqref{twisted} it follows that $\Lam_{L}c\Lam_K=\eps(\Lam_{L})c\Lam_K$ which shows that $L\subset\; ^CK$.
\end{proof}
Note that Remark \ref{nf} together with the previous proposition implies that $\;^CKC\subseteq CK$.
\bc
One has that $\;^CK\subseteq CKC^*$.
\ec
\begin{proof}
Since $S(C)=C^*$ by applying the antipode $S$ to the above inclusion one obtains that $C^*\:^CK\subseteq KC^* $. Therefore $CC^*\;^CK\subseteq CKC^*$ and then one has $\;^CK\subseteq CC^*\;^CK\subseteq CKC^*$.
\end{proof}

\bt\label{larg}
One has that $\;^CK$ is the largest Hopf subalgebra $L$ of $H$ with the property $LC\subseteq CK$.
\et
%
%
%
\begin{proof}
We have seen above that $\;^CKC\subseteq CK$. Suppose now that $LC\subseteq CK$ for some Hopf subalgebra $L$ of $H$. Then by Remark \ref{nf} it follows that $LCK=CK$. Thus  by passing to regular characters one has that $\Lam_Lc\Lam_K=\eps(\Lam_L)c\Lam_K$ which shows the inclusion  $L\subseteq \;^CK$.
\end{proof}
\bp \lb{invcoset} Let $H$ be a semismple Hopf algebra and $K$ be a Hopf subalgebra of $H$.
Then for any subcoalgebra $D$ with $DK=CK$ one has that $\:^DK=\;^CK$.
\ep
\begin{proof}
One has that  $\;^CKCK=CK$. If $D\subset CK$ then by Remark \ref{nf} one has that $\;^CKDK=\;^CKCK=CK=DK$ which shows that $\;^CK=\;^DK$. 
\end{proof}

\section{Mackey type decompositions for representations of Hopf algebras}\lb{mackp}
Let $K$ be a Hopf subalgebra of a semisimple Hopf algebra $H$ and $M$ be a finite dimensional $K$-module.  Note that for any simple subcoalgebra $C$ of $H$ one has by Proposition \ref{modstr} that $\:^CM:=CK\ot_KM$ is a left $\;^CK$-module via the left multiplication with elements of $\;^CK$.
\br
Let $H=kG$ be a group algebra of a finite group $G$ and $K=kA$ for some  subgroup $A$ of $G$. Then note that $\;^CM:=CK\ot_KM$ coincides to the usual conjugate module $\;^gM$ if $C=kg$ for some $g \in G$. Recall that $\;^gM=M$ as vector spaces and $(gag^{-1}).m=a.m$ for all $a \in A$ and all $m \in M$.
\er
\subsection{Proof of Theorem \ref{mack}.}

\bn{proof} Since by definition of the double cosets one has $H=\oplus_{C \in L\backslash H/K}LCK$ and each $LCK$ is a free $K$-module, the following decomposition of $L$-modules follows:
\beq
M\uw^H_K\dw^H_L=H\ot_KM\cong \oplus_{C \in L\backslash H/K}(LCK\ot_KM).
\eeq
Consider now the $k$-linear map $\pi^{(C)}_M: L\ot_{L\cap \;^CK}(CK\ot_KM)\ra LCK\ot_KM$ given by $$l \ot_{L\cap \;^CK}(cx\ot_Km)\mapsto lcx\ot_Km$$ for all $l \in L$, $x \in K$, $c \in C$ and $m \in M$. It is easy to see that $\pi^{(C)}_M$ is a well defined map and clearly a surjective morphism of $L$-modules. Then  $\pi_M:=\oplus_{C \in L\backslash H/K}\pi^{(C)}_M$ is surjective morphism of $L$-modules and the proof is complete.
\end{proof}
\bn{rem}\label{rems}
 Suppose that for $M=k$ one has that $\pi_k$ isomorphism in Theorem \ref{mack}. Then using a dimension argument it follows that the same epimorphism $\pi_M$ from Theorem \ref{mack} is in fact an isomorphism for any finite dimensional left $H$-module $M$.
 \end{rem}
\subsection{Mackey pairs}\lb{mp}
It follows from the proof above that $(L,K)$ is a Mackey pair if and only if $\pi_k$ is an isomorphism, i.e. if and only if  each $\pi^{(C)}_k$ is isomorphism for any simple subcoalgebra $C$ of $H$. Since $\pi^{(C)}_k$ is surjective passing to dimensions one has that $(L,K)$ is a Mackey pair if and only if 
\beq\lb{mtr}
\dim LCK=\frac{(\dim L) \;(\dim CK)}{\dim L \cap \;^CK}
\eeq
for any simple subcoalgebra $C$ of $H$.
\md
Note that for $C=k1$ the above condition can be written as
\beqn
\dim\;LK=\frac{(\dim L)(\dim K)}{\dim(L\cap K)}.
\eeqn
\br
Note also that for any Mackey pair it follows that \beq\frac{\dim \;L \cap \;^DK}{\dim DK}=\frac{ \dim \;L \cap \;^CK}{\dim CK}\eeq if $LCK=LDK$.
\er
\bn{example}
Suppose that $L$, $K$ are Hopf subalgebras of $H$ with $LK=KL$. Then $(L,K)$ is a Mackey pair of Hopf subalgebras of $LK$ by \cite[Proposition 3.3]{coset}.
\end{example} 
\subsection{Proof of Theorem \ref{grpcase}}

Let $L,K$ be two Hopf subalgebras of a semisimple  Hopf algebra $H$ and let $C$ be a simple subcoalgebra of $H$. Note that equations \eqref{twisted} and \eqref{onesided} implies that $LCK$ can be written as a direct sum of right $K$-cosets,
\beq\lb{dc}
LCK=\bigoplus_{DK \in \mtc{S}}DK,
\eeq
for a subset $\mtc{S}\subset H/K$ of right cosets of $K$ inside  $H$. Note that always one has $CK\in \mtc{S}$.
\md
Next we give a proof for the main result Theorem \ref{grpcase}.
\begin{proof}
Suppose that $L=kG$. By Equation \eqref{mtr} one has to verify 
\beq\lb{mtr'}
\dim (kG)CK=\frac{|G| \;(\dim CK)}{\dim kG \cap \;^CK}
\eeq
 for any subcoalgebra $C$ of $H$. Since $kG \cap \;^CK$ is a Hopf subalgebra of $kG$ it follows that $kG \cap \;^CK=kG_C$ for some subgroup $G_C$ of $G$. By Equation \eqref{chcj}  it follows that
$
G_C=\{g \in G\;|\;gd\Lam_K=d\Lam_K\}
$
where $d \in \Irr(H^*)$ is the character associated to $C$. In  terms of subcoalgebras this can be written as 
$
G_C=\{g \in G\;|\;gCK=CK\}.
$

With the above notations Equation \eqref{mtr'} becomes
\beq\lb{mtr''}
\dim (kG)CK=\frac{|G| }{|G_C|}\dim CK
\eeq
 Note that the group $G$ acts transitively on the set $\mtc{S}$ from  Equation \eqref{dc}. The action is given by $g.DK=gDK$ for any $g \in G$ and any $DK\in \mtc{S}$. Let $St_C$ be the stabilizer of the right coset $CK$. Thus the subgroup $St_C$ of $G$ is defined by 
$
St_C=\{g \in G\;|\;gCK=CK\}
$
which shows that $St_C=G_C$. Note that $\dim DK= \dim CK$ for any $DK\in \mtc{S}$ since $DK=gCK$ for some $g \in G$. Thus $\dim (kG)CK = |\mtc{S}|(\dim CK)$ and Equation \eqref{mtr''} becomes 
\beq
|\mtc{S}|=\frac{|G|}{|G_C|}
\eeq
which is the same as the formula for the size of the orbit $\mtc{S}$ of $CK$ under the action of the finite group $G$.
\end{proof}
\section{New examples of Green functors}\lb{newgf}
In this section we construct new examples of Green functors arising from gradings on the category of corepresentations of  semisimple Hopf algebras.

\subsection{Gradings of fusion categories} In this subsection we recall few basic results on gradings of fusion categories from \cite{NG} that will be further used in the paper.
For an introduction to fusion categories one might consult \cite{ENO}. Let $\cc$ be a fusion category and $\co(\mtc{C})$ be the set of isomorphism classes of simple objects of $\cc$. Recall that the fusion category $\cc$ is graded by a finite group $G$ if there is a function $\deg: \co(\cc) \ra G$ such that for any two simple objects $X,Y \in \co(\cc)$ one has that $\deg(Z)=\deg(X)\deg(Y)$ whenever $Z\in \co(\cc)$ is a simple object such that $Z$ is a constituent of $X\ot Y$.
Alternatively, there is a decomposition $\mtc{C}=\oplus_{g \in G}\mtc{C}_g$
such that the the tensor functor of $\cc$ sends $\cc_g \ot \cc_h$ into $\cc_{gh}$. Here $\cc_g$ is defined as the full abelian subcategory of $\cc$ generated by the simple objects $X$ of $\cc$ satisfying $\deg(X)=g$. Recall that a grading is called universal if any other grading of $\cc$ is arising as a quotient of the universal grading. The universal grading always exists and its grading group denoted by $U_{\cc}$ is called the universal grading group.
\br\lb{hopf}
If $\cc=\rep(H)$  for a semisimple Hopf algebra $H$ then by \cite[Theorem 3.8]{NG} it follows that the Hopf center (i.e. the largest central Hopf subalgebra) of $H$ is $kG^*$ where $G$ is the universal grading group of $\cc$.
We denote this Hopf center by $\mtc{HZ}(H)$. Therefore one has  $\mtc{HZ}(H)=kG^*$ where $G=U_{\rep(H)}$. 
Moreover, in this case, by the universal property any other grading on $\cc=\rep(H)$ is given  by a quotient group $G/N$ of $G$. The corresponding graded components of $\cc$ are given by 
\beq
\cc_{\bar{g}}=\{M\in \Irr(H)\;|\;M\dw^H_{k^{G/N}}=(\dim M) \bar{g}\}
\eeq
for all $g \in G$. Here $k^{G/N}\subset k^{G}$ is regarded as a central Hopf subalgebra of $H$. Also note that in this situation one has a central extension of Hopf algebras:
\bq
k\ra k^{G/N}\ra H \ra H//k^{G/N}\ra k.
\eq
\er
\subsection{Gradings on $\rep(H^*)$ and cocentral extensions}
Suppose that $H$ is a semisimple Hopf algebra such that the fusion category $\rep(H^*)$ is graded by a finite group $G$. Then the dual version of Remark \ref{hopf} implies that $H$ fits into a cocentral extension
\beq\lb{coc}
k \ra B \ra H \xra{\pi} kG \ra k.
\eeq 
Recall that such an exact sequence of Hopf algebras is called cocentral if $kG^*\subset \mtc{Z}(H^*)$ via the dual map $\pi^*$. On the other hand, using the reconstruction theorem from \cite{AD} it follows that 
\beq\lb{andc}
H \cong B\;^{\tau}\#_{\sg} \;kF
\eeq
for some cocycle $\sg:B\ot B \ra kF$ and some dual cocyle $\tau:kF \ra B\ot B$.
\md
 For  any such cocentral sequence it follows  that $G$ acts on $\rep(B)$ and by \cite[Proposition 3.5] {natalecoc} that $\rep(H)=\rep(B)^G$, the equivariantized fusion category. For the main properties of group actions and equivariantized fusion categories one can consult for example \cite{N}. Recall that the above action of $G$ on $\rep(B)$ is given by $  T:G\xra{} \mtr{Aut}_{\ot}(\rep(B))$, $g \mapsto T^g$. For any $M \in \rep(B)$ one has that $T^g(M)=M$ as vector spaces and the action of $B$ is given by
$b.^gm:=(g.b).m$ for all $g \in G$ and all $b \in B$, $m \in M$. Here the weak action of $G$ on $B$ is the  action used in the crossed product from Equation \eqref{andc}.

For any subgroup $M$ of $G$ it is easy to check that $H(M)=B\#_{\sg}kM$, i.e. $H(M)$ is the unique Hopf subalgebra of $H$ fitting the exact cocentral sequence
\beq\lb {ses}
k \ra B \ra H(M) \ra kM \ra k.
\eeq 

\bl\lb{grd*} Let $H$ be a semisimple Hopf algebra. Then gradings on the fusion category $\rep(H^*)$ are in one-to one correspondence with cocentral extensions
\beq\lb{sesi}
k \ra B \ra H \xra{\pi} kG\ra k.
\eeq 
\el
\bpf
We have shown at the beginning of this subsection how to associate a cocentral extension to any $G$-grading on $\rep(H^*)$.

Conversely, suppose one has a cocentral exact sequence as in Equation \eqref{sesi}. 
Then $\rep(H^*)$ is graded by $G$ where the graded component of degree $g\in G$ is given by
\beq
\rep(H^*)_g=\{d \in \Irr(H^*)\;|\;\pi(d)=\eps(d)g\}.
\eeq

Indeed, since $k^G\subset \mtc{Z}(H^*)$ via $\pi^*$ it follows that $k^G$ acts by scalars on each irreducible representation of $H^*$. Therefore for any $d \in \Irr(H^*)$ one has $d\dw^{H^*}_{k^G}=\eps(d)g$ for some $g \in G$. It follows then by  \cite [Theorem 3.8]{NG})  that $\rep(H^*)$ is $G$-graded and 
\beq
\rep(H^*)_g=\{d \in \Irr(H^*)\;|\;d\dw^{H^*}_{k^G}=\eps(d)g\}
\eeq
On the other hand it is easy to check that one has $\pi(d)=d\dw^{H^*}_{k^G}$ for any $d \in \Irr(H^*)$ (see also Remark 3.2 of \cite{coc}.)

Clearly the two constructions are inverse one to the other.
\epf
\subsection{New examples of Mackey pairs of Hopf subalgebras}
Let $H$ be a semisimple Hopf algebra and $\cc=\rep(H^*)$. Since $H^*$ is also a semisimple Hopf algebra \cite{Lard} it follows that $\cc$ is a fusion category. For the rest of this section fix an arbitrary $G$-grading $\cc=\oplus_{g \in G}\cc_g$ on $\cc$.

For any subset $M\subset G$ define $\cc_M:=\oplus_{m \
\in M}\cc_m$ as a full abelian subcategory of $\cc$. Thus $\co(\cc_M)=\sqcup_{m \in M}\co(\cc_m)$. Let also $H(M)$ to be the subcolagebra of $H$ generated by 
all the simple subcoalgebras of $H$ whose irreducible $H^*$-characters belong to $\co(\cc_M)$.

For any subcoalgebra $C$ of $H$ denote by $\Irr(C^*)$ the irreducible characters of the dual algebra $C^*$. Therefore by its definition $H(M)$ verifies the equality $\Irr(H(M)^*)=\co(\cc_M)$ and as a coalgebra can be written as
$H(M)=\bigoplus_{\{d \in \co(C_m)\;|\;m\in M\}}C_d.$ Note that if $M$ is a subgroup of $G$ then $H(M)$ is a Hopf subalgebra of $H$ by Remark \ref{princ}.
\md
For any simple subcoalgebra $C$ of $H$ whose associated irreducible character $d \in \Irr(H^*)$ has degree $g$ we will also write for shortness that $\deg(C)=g$.
\bp \lb{dcc} Let $H$ be semisimple Hopf algebra and $G$ be the universal grading group of $\rep(H^*)$.
Then for any arbitrary two subgroups $M$ and $N$ of $G$, the set of double cosets $H(M)\backslash H \slash H(N)$ is canonically bijective to the set of group double cosets $M\backslash N\slash G$. Moreover, the bijection is given by $H(M)CH(N)\mapsto M\deg(C)N$.
\ep
\begin{proof} By Remark \ref{princ} one has the following equality in terms of irreducible $H^*$-characters: 
\beqn 
\Irr (H(M)CH(N)^*)=\co(\cc_{M\deg(C)N})
\eeqn 
Thus if $H(M)CH(N)=H(M)DH(N)$ then $\deg(C)=\deg(D)$ which shows that the above map is well defined. Clearly the map $H(M)CH(N)\mapsto M\deg(C)N$ is surjective. The injectivity of this map also follows from Remark \ref{princ}.
\end{proof}

Note that the proof of the previous Proposition implies that the coset $H_x=H(M)CH(N)$ with $\deg(C)=x$ is given by 
\beq\lb{df2}
H_x=\bigoplus_{\{d \in \co(\cc_{mxn})\;|\; m \in M, \; n\in N\}}C_d.
\eeq 

\bp\lb{eqn}
Suppose that $V \in H(M)-mod$, i.e. $V$ is a $B\#_{\sg}kM$-mdule. Then as $B$-modules one has that 
$\;^CV\cong T^{g^{-1}}(\mtr{Res}^{H(M)}_B(V))$
where $g \in G$ is chosen such that $\deg(C)=g$. Moreover $\;^C(V\ot W)\cong  \;^CV\ot \;^CW$ for any two left $H(M)$-modules $V$ and $W$.
\ep
\ncm{\dc}{\#_{\sg}}
\bpf
Note that in this situation one has that $\;^CH(M)=H(\;^gM)=B\#_{\sg}k^gM$. By definition one has $\;^CV=CH(M)\ot_{H(M)}V=H(gM)\ot_{H(M)}V$. Thus
\beq
\;^CV=(B\# kgM)\ot_{B \#kM}V\cong  kgM\ot_{kM}V
\eeq
where the inverse of the last isomorphism is given by $g\ot_{kM}v \mapsto (1\#g)\ot_{ B \#kM}v$. Note that $B$ acts on $kgM\ot_{kM}V$ via $b.(g\ot_{kM} v)=g\ot_{kM} (g\inv.b)m$  for all $b \in B$, $v\in V$. This shows that indeed $\;^CV\cong T^{g\inv}(\mtr{Res}^{H(M)}_B(V))$ as $B$-modules. Moreover it follows that $\;^CV$ can be identified to $V$ as vector spaces with the $B\#_{\sg}kgMg\inv$-module structure given by 
$b.v=(g\inv.b)v$ and  
$
(ghg\inv).v=([g\inv.(\sg(ghg\inv, g)\sg\inv(g,h))]\#_{\sg}h).v.
$
for all $g \in G$, $h \in M$ and $v \in V$. Then it can be checked by direct computation that the map $v\ot w \mapsto \tau^{-1}(g)(v\ot w)$ from \cite[Proposition 3.5] {natalecoc} is in this case a morphism of $B\#_{\sg}kgMg\inv$-modules. In order to do that one has to use the compatibility conditions from Theorem 2.20 of \cite {AD}.

\epf
\subsection{Examples of Mackey pairs arising from group  gradings on the category  $\rep(H^*)$.}
Let as above $H$ be a semisimple Hopf algebra with
$\cc=\os_{g \in G}\cc_g$
be a group grading of $\cc:=\rep(H^*)$. It follows that 
\beq\lb{fpd}
\fp(\cc_g)=\frac{\dim H^*}{\dim \mtc{HZ}(H^*)}
\eeq 
for all $g \in G$ where $\fp(\cc_g):=\sum_{V \in \co(\cc_G)}(\dim V)^2$ is the Perron-Frobenius dimension of the full abelian subcategory $\cc_g$ of $\cc$. 

\bt\lb{grd} Let $H$ be a semisimple Hopf algebra and $M$, $N$ be any two subgroups of $G$. With the above notations the pair $(H(M),H(N))$ is a Mackey pair of Hopf subalgebras of $H$.
\et

\begin{proof} Put $L:=H(M)$ and $K:=H(N)$. Therefore $\Irr(L^*)=\co(\cc(M))$ and $\Irr(K^*)=\co(\cc(N))$. Then we have to verify Equation \eqref{mtr} for any simple subcoalgebra $C$. Fix a simple subcoalgebra $C$ of $H$ with $\deg(C)=x$.
As above one has $\;^CH(M)=H(\;^xM)$.

It is easy to verify that $H(M)\cap H(N)=H(M\cap N)$ for any two subgroups $M$ and $N$ of $G$. This implies that $L\cap \;^CK=H(N\cap\;^xM)$. \onh\;from Equation \eqref{df2} note that $\dim LCK=|MxN|\fp(\cc_1)$.

Then Equation \eqref{mtr} is equivalent to the well known formula for the size of a double coset relative to two subgroups: 
\beq
|MxN|=\frac{|M||N|}{|M\cap \;^xN|},
\eeq
for any $x \in G$.
\end{proof}
\br The fact that $(H(M), H(N))$ is a Mackey pair also follows in this case from a more general version of Mackey's decomposition theorem that holds for the action of any finite group on a fusion category. This results will be contained in a future paper of the author. 
\er
\br
It also should be noticed that the author is not aware  of any pair of Hopf subalgebras that is not a Mackey pair. It would be interesting to construct such counterexamples if they exist.
\er
\subsection{Mackey and Green functors}For a finite group $G$ denote by $\cs(G)$ the lattice of all subgroups of $G$. Following \cite{thv} a Mackey functor for $G$ over a ring $R$ can be regarded as a collection of vector spaces $M(H)$ for any $H \subset \cs(G)$ together with a family of morphisms
$I^{L}_{K} : M(K)\ra M(L)$,
$R^{L}_{K}: M(L)\ra M(K)$, and
$c_{K, g} : M(K)\ra M( \;^{g}K)$
for all subgroups $K$ and $L$ of $G$ with $K\subset L$ and for all  $g \in G$. This family of morphisms has to satisfy the following compatibility conditions:
\bne\item\lb{m0}
$I_{H}^H, \;R_{H}^H, c_{H, h} : M(H)\ra M(H)$ are the identity morphisms for all subgroups $H$ of $G$ and any $h \in H$,
\item \lb{m1}
$R_{K}^{J}R_{H}^{K}=R^{J}_{K}$, for all subgroups $J\subset K\subset H$,
\item \lb{m2}
$I_{H}^{K}I_{J}^{H}=I^{K}_{J},
$
for all subgroups $J\subset K\subset H$,
\item\lb{m3} 
$c_{K, g}c_{K, h} = c_{K, gh} 
$
for all elements $g, h \in G$.
\item \lb{maci}For any three subgroups $J, L \subseteq K$ of $G$ and any $a \in M(J)$ one has the following Mackey axiom:
\beqn
R^J_L(I_J^K(a))=\sum_{x \in J\backslash K\slash L}I^L_{L\cap \;^xJ}(R^{\;^xJ}_{\;^xJ\cap L}(c_{J, x}(a)))
\eeqn

\hskip -1cm Moreover, a Green functor is a Mackey functor $M$
such that for any subgroup $K$ of $G$ one has 

\hskip -1cm that $M(K)$
is an associative $R$-algebra with identity and the following conditions are satisfied:
\item \lb{g1} $R_{K}^{L}\;\;
\text{and} \;\; c_{{K, g}}\;\;\text{are always unitary R-algebra homomorphisms,}
$
\item  \lb{g2}
$
I_{K}^{L}(aR^{L}_{K}(b)) = I_{K}^{L}(a)b,
$
\item\lb{g3}
$I_{K}^{L}(R^{L}_{K}(b)a) =b I_{K}^{L}(a)
$
 for all subgroups $K\subseteq L\subseteq G$ and all $a \in M(K)$ and $b \in M(L)$. 
 \ene
 \md
Green functors play an important role in the representation theory of finite groups (see for example \cite{thv}).
\subsection{New examples of Green functors}
Next Theorem allows us to construct new examples of Green functors from semisimple Hopf algebras.
\bt
Let $H$ be a semisimple Hopf algebra and $G$ be a grading group for the fusion category $\rep(H^*)$. Then the functor $M \mapsto K_0(H(M))$ is a Green functor.
\et
\ncm{\modd}{\mtr{mod}}
\bpf
By Proposition \ref{grd*} there is a cocentral extension
\beq\lb{sesi}
k \ra B \ra H \xra{\pi} kG\ra k
\eeq 
for some Hopf subalgebra $B\subset H$. Then as above, for a simple subcoalgebra $C$ of $H$ with associated character $d \in H^*$ one has that if $\pi(d)=g$ for some $g \in G$ then $\pi(C)=kg$.

The map $R^L_K: K_0(H(L)) \ra K_0(H(K))$ is induced by the restriction map $\mtr{Res}^{H(L)}_{H(K)}:H(L)-\modd\ra H(K)-\modd$. Similarly, the map $I^L_K$ is induced by the induction functor between the same two categories of modules.  Clearly $R^L_K$ is a unital algebra map and the compatibility conditions \ref{g2} and \ref{g3} follow from the adjunction of the two functors. Moreover conditions \ref{m1} and \ref{m2} are automatically satisfied.

Define $c_{L,g}:K_0(L)\ra K_0(\;^gL)$ by $[M]\mapsto [\;^CM]$ where $C$ is any simple subcoalgebra of $H$ chosen with the property that $\deg(C)=g$. It follows by Proposition \ref{eqn}  that $c_{L, g}$ is a well defined algebra map. Condition \ref{m3} is equivalent to $T^{gh}(M)\cong T^{g}T^h(M)$ which is automatically satisfied for a group action on a fusion category.

It is easy to see that all other axioms from the definition of a Green functor are satisfied. For example, the Mackey decomposition axiom \ref{maci} is satisfied by Theorem \ref{grd}.
\epf
\section{On normal Hopf subalgebras of semisimple Hopf algebras}\lb{nrm}Recall that a Hopf subalgebra $L$ of a Hopf algebra $H$ is called a normal Hopf subalgebra if it is stable under the left and right adjoint action of $H$ on itself. When $H$ is a semisimple Hopf algebra, since $S^2=id$, in order for $L$ to be normal,  it is enough to be closed only under the left adjoint action, i.e. $h_1LS(h_2) \subset L$ for any $h \in H$. Let $L^+:=L \cap \ker\;\eps$ and set $H//L:=H/HL^+$. Since $HL^+$ is a Hopf ideal of $H$ (see for example \cite{montg}) it follows that $H//L$ is a quotient Hopf algebra of $H$. Moreover $(H//L)^*$ can be regarded as a Hopf subalgebra of $H^*$ via the dual map of the canonical Hopf projection $\pi_L:H\ra H//L$.
\bp\lb{norm}
Suppose that $H$ is a semisimple Hopf algebra. Then for any normal Hopf subalgebra $K$ of $H$ one has that $(K,K)$ is a Mackey pair of Hopf algebras.
\ep
\begin{proof}
Note $KC=CK$ for any subcoalgebra $C$ of $K$ since $K$ is a normal Hopf subalgebra of $H$.  Then for any simple subcoalgebra $C$ of $H$ the dimension condition from Equation \eqref{mtr} can be written as
\beq\lb{mtrn}
\dim CK=\frac{(\dim K) \;(\dim CK)}{\dim K \cap \;^CK}
\eeq
which is equivalent to $K \cap \;^CK=K$. This equality follows by the third item of Remark \eqref{comm}.
\end{proof}
\subsection{Irreducibility criterion for an induced module}
\br Let $G$ be a finite group and $H$ be a normal subgroup $H$ of $G$
Then \cite[Corollary 7.1]{Se76} implies that an induced module $M\uw^G_H$ is irreducible if and only if $M$ is irreducible and $M$ is not isomorphic to any of its conjugate module $\;^gM$.
\er
Previous Theorem allows us to prove the following Proposition which is an improvement of Proposition 5.12 of \cite{coset}. The second item is also a generalization of \cite[Corrolary 7.1]{Se76}.
\bp\lb{irredind}
Let $K$ be a normal Hopf subalgebra of a \sem Hopf algebra $H$ and $M$ be a finite dimensional $K$-module. \bne
\item Then \beqn
M\uw^H_K\dw^H_K\cong \bigoplus_{C\in H/ K} \;^CM
\eeqn
as $K$-modules.
\item $M\uw^H_K$ is irreducible if and only if $M$ is an irreducible $K$-module which is not a direct summand of any conjugate module $\;^CM$ for any simple subcoalgebra $C$ of $H$ with $C\not\subset K$.
\ene
\ep
\begin{proof}
\bne
\item Previous Proposition implies that 
\beq
M\uw^H_K\dw^H_K\cong \bigoplus_{C\in K\backslash H/K} K\ot_{K\cap \;^CK}\;^CM
\eeq
as $K$-modules. On the other hand since $K$ is normal note that $CK=KC$ and therefore the space $K\backslash H/K$ of double cosets coincides to the space $H/ K$ of left (right) cosets (see also Subsection \ref{not} for the notation). In the the proof of the same Proposition \ref{norm} it was also remarked that $K\cap \;^CK=K$.
\item
One has that $M\uw^H_K$ is an irreducible $H$-module if and only if $$\dim_k \mtr{Hom}_H(M\uw^H_K,\;M\uw^H_K)=1.$$ Note that by the Frobenius reciprocity one has the following $\mtr{Hom}_H(M\uw^H_K,\;M\uw^H_K)=\mtr{Hom}_K(M,\;M\uw^H_K\dw^H_K)$. Then previous item implies that 
\beq
\mtr{Hom}_K(M,\;M\uw^H_K)\cong \bigoplus_{C\in H/K} \mtr{Hom}_K(M, \;^CM)
\eeq
Since for $C=k$ one has $\;^kM=M$ it follows that $\mtr{Hom}_K(M, \;^CM)=0$ for all $C\not\subset K$.
\ene
\end{proof}
\subsection{A tensor product formula for induced representations}\lb{tpf}
We need the following preliminary tensor product formula for induced representations which appeared in \cite{bkk}.
\bp\lb{kul}  Let $K$ be a Hopf subalgebra of a semisimple Hopf algebra $H$. Then for any $K$-module $M$ and any $H$-module $V$ one has that 
\beq
M \uw^H_K\otimes V\cong (M \otimes V\dw^{\;H}_{\;K})\uw^{\;H}_{\;K}
\eeq
\ep

{\bf Proof of  Theorem \ref{tp}:}
Applying Proposition \ref{kul} one has that
\beq
M\uw_K^H\ot N\uw_L^H  \cong (M\uw_K^H\dw_L^H\ot N) \uw_L^H 
\eeq
On the other hand, by Theorem \ref{mack} one has 
\beq
M\uw_K^H\dw_L^H\cong  \oplus_{C \in L\backslash H/K}(L\ot_{L\cap \;^CK}(CK\ot_KM))
\eeq
Thus, 
\bn{eqnarray*}
M\uw_K^H\ot N\uw_L^H  & \cong &   (M\uw_K^H\dw_L^H\ot N) \uw_L^H \\ &\xrightarrow{\cong} &  \oplus_{C \in L\backslash H/K}((L\ot_{L\cap \;^CK}(CK\ot_KM))\ot N)\uw_L^H\end{eqnarray*}

Applying again Proposition \ref{tp} for the second tensor product one obtains that
\bn{eqnarray*}  M\uw_K^H\ot N\uw_L^H &\xrightarrow{\cong} & ((CK\ot_KM)\ot N\dw^{\;L}_{\;L\cap \;^CK})\uw_{\;L\cap \;^CK}^L\uw^H_L \\ &\xrightarrow{\cong} & \oplus_{C \in L\backslash H/K}H\ot_{L\cap \;^CK}((CK\ot_KM)\ot N\dw^{\;L}_{\;L\cap \;^CK}).
\end{eqnarray*}
\br
Note that the above Theorem always applies for $K=L$ a normal Hopf subalgebra of $H$.
\er
\bibliographystyle{amsplain}
\bibliography{v1}
\address{
   \email{sebastian.burciu@imar.ro}}
\end{document}